\documentclass{tac}

 \usepackage{hyperref}
\usepackage {xspace}
 \input diagxy
 \title{Topological $*$-autonomous categories, revisited}
 \author{Michael Barr}
 \address{McGill University, Dept. of Math. and Stats.\\ 805 Sherbrooke
W, Montreal, QC H3A 2K6}
 \eaddress{mbarr@barrs.org}

\thanks{I would like to thank the NSERC of Canada}
\keywords{duality, Chu construction, Mackey spaces}
\amsclass{18D15, 22D35, 46A20}

\let\cont\supseteq
\def\nbd{neighborhood\xspace}
\def\incto{\to/into->/}
\mathrmdef{id}

\def\makeother#1{\catcode`#1=12}
{\catcode`\^^M=12 \endlinechar=-1 %
   \gdef\xcomment#1^^M{\def\test{#1}
       \ifx\test\endcomment \let\next=\endgroup
       \else \let\next=\xcomment \fi
    \next}
}
\def\comment{\begingroup
   \let\do=\makeother \dospecials
   \endlinechar`\^^M \catcode`\^^M=12 \xcomment }
{\escapechar=-1 \xdef\endcomment{\string\\endcomment}
}

\mathrmdef{Vect}
\mathzcdef[Asc]A
\mathzcdef[Bsc]B
\mathzcdef[Csc]C
\mathzcdef[Dsc]D
\mathzcdef[Esc]E
\mathzcdef[Msc]M
\mathzcdef[Ssc]S
\mathzcdef[Tsc]T
\mathzcdef[Vsc]V
 \def\inv{^{-1}}
\let\eps\epsilon
\let\lam\lambda
\def\|#1{|\!|#1|\!|}
\DeclareFontFamily{OT1}{pzc}{}
\DeclareFontShape{OT1}{pzc}{m}{it}%
     {<-> s * [1.100] pzcmi7t}{}
\DeclareMathAlphabet{\mathscr}{OT1}{pzc}{m}{it}
\DeclareSymbolFont{script}{OT1}{pzc}{m}{it}
\DeclareFontFamily{OT1}{cmssbx}{}
\DeclareFontShape{OT1}{cmssbx}{m}{n}{
            <5> <6> <7> <8> 
   <9> <10> <10.95> <12> <14.4> <17.28> <20.74> <24.88> cmssbx10 }{}

\DeclareSymbolFont{ssbx}{OT1}{cmss}{bx}{n}
\def\inv{^{-1}}
\let\thm\theorem
\let\prop\proposition
\let\cor\corollary
\let\lem\lemma
\let\eth\endtheorem
\let\pf\proof
\let\epf\endproof

\let\defn\definition

\def\o{\mathop{\raise.3ex\hbox{$\scriptscriptstyle\circ$}}}
\let\x=\times
\let\inc\subseteq
\let\cont\supseteq
\let\iso=\cong
\let\cong=\equiv

\mathchardef\gt="313E 
\mathchardef\lt="313C 

\def\cahiers#1{Cah\-iers de Topo\-lo\-gie et G\'e\-o\-m\'e\-trie
Dif\-f\'e\-ren\-ti\-elle\ifnum#1>25 { Cat\-\'e\-go\-rique}\fi, {\bf #1}}
\def\tac#1{Theory and Applications of Categories, {\bf #1}}

\def\mscs#1{Mathematical Structures Computer Science, {\bf #1}}

 \mathopdef{Hom}
 \mathopdef{Mod}
 \mathrmdef{chu}
 \mathrmdef{Chu}
 \mathrmdef{Id}
 \mathssbxdef{R}
 \mathssbxdef{N}
 \mathssbxdef{C}
 \mathssbxdef{Z}
 \mathssbxdef{Q}
 \mathssbxdef{S}
 \mathzcdef{Ab}
\def\adj{\mathrel{\raise.78ex\hbox{\vrule width.9em height.2pt depth.2pt
  \vrule width.4pt height.55ex depth.55ex}}}
\mathcode`\<="4268 
\mathcode`\>="5269 
\let\+\oplus
\let\*\otimes

 \let\phi\varphi
 \def\nbd{neighbourhood\xspace}
 
\def\r{\mathchoice
{\mathinner{\relbar\mkern-7.5mu\circ}}
{\mathinner{\relbar\mkern-7.5mu\circ}}
{\mathinner{\scriptstyle\relbar\mkern-3.5mu\circ}}
{\mathinner{\scriptscriptstyle\relbar\mkern-3.5mu\circ}}
}
\begin{document}

\maketitle

\abstract Given an additive equational category with a closed symmetric
monoidal structure and a potential dualizing object, we find sufficient
conditions that the category of topological objects over that category
has a good notion of full subcategories of strong and weakly topologized
objects and show that each is equivalent to the chu category of the
original category with respect to the dualizing object.\endabstract

\section{Introduction}

This paper is an updated version of \cite{top}.  In the course of
preparing some lectures on the subject, I discovered to my great chagrin
that that paper was badly flawed.  Several of the arguments had gaps or
flaws.  In the process of repairing them, I discovered that the main
results were not only correct, but that better results were
available. This updated paper is the result.

\cite{M} introduced the category of pairs of vector spaces,
equipped with a bilinear pairing into the ground field.  It is likely
that he viewed this abstract duality as a replacement for the topology.
See also \cite{M46}, the review of the latter paper by Dieudonn\'e as
well as Dieudonn\'e's review of \cite{A}, for a clear expression of this
point of view.  In this paper we fully answer this question.

\cite{B00} showed that the full subcategory of the category of
(real or complex) topological vector spaces that consists of the Mackey
spaces (defined in \ref{strong&weak} below)is $*$-autonomous and
equivalent to both the full subcategory of
weakly topologized topological vector spaces and to the full subcategory
of topological vector spaces topologized with the strong, or Mackey
topology.  This means, first, that those subcategories can, in
principle at least, be studied without taking the topology into
consideration.  Second it implies that both of those categories are
$*$-autonomous.

Andr\'e Joyal has raised the question whether there was a similar result
for vector spaces over the field $\Q_p$ of $p$-adic rationals.  This was
mentioned in the original paper, but oddly that situation was not
explored.  Thinking about this question, I realized that there is a
useful general theorem that answers this question for any locally
compact field and also for locally compact abelian groups.

The results in this paper support the following conclusion.  Let $K$ be
a {\bf spherically complete field} (defined below) and $|K|$ its
underlying discrete field.  Then the following five categories are
equivalent:

\begin{enumerate}

\item$\chu(K$-$\Vect,|K|)$ (Section \ref{chu})

\item The category $\Vsc_w(K)$ of topological $K$-spaces topologized
with the weak topology for all their continuous linear functionals into
$K$.

\item The category $\Vsc_s(K)$ of topological $K$-spaces topologized
with the strong topology (see Section \ref{weakstrong}) for all their
continuous linear functionals into $K$.

\item The category $\Vsc_w(|K|)$ of topological $|K|$-spaces topologized
with the weak topology for all their continuous linear functionals into
$|K|$.

\item The category $\Vsc_s(|K|)$ of topological $|K|$-spaces topologized
with the strong topology (see Section \ref{weakstrong}) for all their
continuous linear functionals into $|K|$.
\end{enumerate}

A normed field is {\bf spherically complete} if any family of closed
balls with the finite intersection property has non-empty intersection.
A locally compact field is spherically complete (so this answers Joyal's
question since $\Q_p$, along with its finite extensions, is locally
compact) and spherically complete is known to be strictly strong than
complete.

\cite[Section 2]{top} is a result on adjoint functors that is
interesting and possibly new.  The argument given there is flawed.
Although the result is bypassed in the current paper, the result seemed
interesting enough to give a full proof of.  This appears as an appendix
to this paper.

\subsection{Terminology}
 We assume that all topological objects are Hausdorff.  As we will see,
each of the categories contains an object $K$ with special properties.
It will be convenient to call a morphism $V\to K$ a {\bf functional} on
$V$.  In the case of abelian groups, the word ``character'' would be
more appropriate, but it is convenient to have one word.  In a similar
vein, we may refer to a mapping of topological abelian groups as
``linear'' to mean additive.  We will be dealing with topological
objects in categories of topological vector spaces and abelian groups.
If $V$ is such an object, we will denote by $|V|$ the underlying vector
space or group.

If $K$ is a topological field, we will say that a vector space is
{\bf linearly discrete} if it is a categorical sum of copies of the
field.

\section{The strong and weak topologies}\label{weakstrong}

\subsection{Blanket assumptions.}\label{blanket} Throughout this
section, we make the following assumptions.

\begin{enumerate}

\item $\Asc$ is an additive equational closed
symmetric monoidal category and $\Tsc$ is the category of topological
$\Asc$-algebras.

\item $K$ is a uniformly complete object of $\Tsc$.

\item there is a \nbd $U$ of 0 in $K$ such that \begin{enumerate}{}{}

 \item $U$ contains no non-zero subobject;

 \item whenever $\phi:T\to K$ is such that $\phi\inv(U)$ is
open, then $\phi$ is continuous.
 \end{enumerate}
 \end{enumerate}

\lem \label{main}Suppose there is an embedding $T\incto\prod_{i\in
I}T_i$ and there is a morphism $\phi:T\to K$.  Then there is a finite
subset $J\inc I$ and a commutative diagram
 $$\bfig
 \Dtriangle<400,500>[T`T_0`K;\phi``\phi_0]
 \morphism(0,1000)/into->/<1000,0>[T`\prod_{i\in I}T_i;]
 \morphism(400,500)/into->/<600,0>[T_0`\prod_{j\in J}T_j;]
 \morphism(1000,1000)<0,-500>[\prod_{i\in I}T_i`\prod_{j\in J}T_j;]
 \efig$$
 Moreover, we can take $T_0$ closed in $\prod_{j\in J}T_j$.\eth

\pf Since $\phi\inv(U)$ is a \nbd of 0 in $T$, it must be the meet with
$T$ of a \nbd of 0 in $\prod_{i\in I}T_i$.  From the definition of the
product topology, we must have a finite subset $J\inc I$ and \nbd{s}
$U_j$ of 0 in $ T_j$ such that
 $$\phi\inv(U)\cont T\cap(\prod_{j\in J}U_j\x\prod_{i\in I-J}T_i)$$
 It follows that
 $$U\cont\phi(T\cap(\prod_{j\in J}0\x\prod_{i\in I-J}T_j))$$
 But the latter is a subobject and therefore must be 0.  Now let
 $$T_0={T\over T\cap(\prod_{j\in J}0\x\prod_{i\in I-J}T_j)}$$
 topologized as a subspace of $\prod_{j\in J}T_0$ and $\phi_0$ be the
induced map.  It is immediate that $\phi_0\inv(U)\cont\prod_{j\in J}U_j$
which is a \nbd of 0 in the induced topology and hence $\phi_0$ is
continuous.  Finally, since $K$ is complete, we can replace $T_0$ by its
closure in $\prod_{j\in J}T_j$.\epf

\thm \label{inj}Suppose $\Ssc$ is a full subcategory of $\Tsc$ that is
closed under finite products and closed subobjects and that $K\in \Ssc$
satisfies the assumptions in \ref{blanket}.  If $\Vsc$ is the closure of
$\Ssc$ under all products and all subobjects and $K$ is injective in
$\Ssc$, then it is also injective in $\Vsc$.\eth

\pf It is sufficient to show that if $V\inc \prod_{i\in I}S_i$ with each
$S_i\in\Ssc$, then every morphism $V\to K$ extends to the product.  But
the object $V_0$ constructed in the preceding lemma is a closed suboject
of $\prod_{j\in J}S_j$ so that $V_0\in\Ssc$ and the fact that $K$ is
injective in $\Ssc$ completes the proof.\epf

\defn A bijective morphism $V\to V'$ in $\Vsc$ is called a {\bf weak
isomorphism} if the induced $\Hom(V',K)\to\Hom(V,K)$ is a bijection.\eth

Of course, a bijective morphism induces an injection so the only issue
is whether the induced map is a surjection.

\prop\label{finprod} A finite product of weak isomorphisms is a weak
isomorphism.\eth

\pf Assume that $J$ is a finite set and for each $j\in J$, $V_j\to V_j'$
is a weak isomorphism.  Then since finite products are the same as
finite sums in an additive category, we have
 $$\mldepf{}&\Hom(\prod V_j',K)\iso\Hom(\sum V_j',K) \iso\prod
\Hom(V_j',K)\\ \iso\prod\Hom(V_j,K)\iso \Hom(\sum V_j,K)\iso \Hom(\prod
V_j,k ) $$

\thm\label{s&t} Assume the conditions of Theorem \ref{inj} and also
suppose that for every object of $\Ssc$, and therefore of $\Vsc$, there
are enough functionals to separate points.  Then for every object $V$ of
$\Vsc$, there are weak isomorphisms $\tau V\to V\to\sigma V$
with the property that $\sigma V$ has the coarsest topology that has the
same functionals as $V$ and $\tau V$ has the finest topology that has
same functionals as $V$.\eth

\pf The argument for $\sigma$ is standard.  Simply retopologize $V$ as a
subspace of $K^{\Hom(V,K)}$.

Let $\{V_i\to V\}$ range over the isomorphism classes of weak
isomorphisms to $V$. We define
$\tau V$ as the pullback in
 $$\bfig
 \Square[\tau V`\prod V_i`V`V^I;```]
 \efig$$
 The bottom map is the diagonal and is a topological embedding so that
the top map is also a topological embedding.  We must show that every
functional on $\tau V$ is continuous on $V$.  Let $\phi$ be a functional
on $\tau V$.  From injectivity, it extends to a functional $\psi$ on
$\prod V_i$.  By Lemma \ref{main}, there is a finite subset $J\inc I$
and a functional $\psi_0$ on $\prod_{j\in J} V_j$ such that $\psi$ is
the composite $\prod_{i\in I}V_i\to\prod_{j\in J}V_j\to^{\psi_0}K$.
Thus we have the commutative diagram
 $$\bfig
 \node 1a(0,600)[\tau V]
 \node 1b(600,600)[\prod_{i\in I}V_i]
 \node 1c(1300,600)[\prod_{j\in I}V_j]
 \node 2(1800,300)[K]
 \node 3a(0,0)[V]
 \node 3b(600,0)[V^I]
 \node 3c(1300,0)[V^J]
 \arrow[1a`1b;]
 \arrow[1b`1c;]
 \arrow[1c`2;\psi_0]
 \arrow[1a`3a;]
 \arrow[1b`3b;]
 \arrow[1c`3c;]
 \arrow[3a`3b;]
 \arrow[3b`3c;]
 \arrow/-->/[3c`2;]
 \efig$$
 The dashed arrow exists because of Proposition \ref{finprod}, which
completes the proof.\epf

\remark\label{strong&weak}
We will call the topologies on $\sigma V$ and $\tau V$ the {\bf weak}
and {\bf strong} topologies, respectively.  They are the coarsest and
finest topology that have the same underlying $\Asc$ structure and the
same functionals as $V$.  The strong topology is also called the {\bf
Mackey topology}.\eth

\prop Weak isomorphisms are stable under pullback.\eth

\pf  Suppose that
 $$\bfig
 \Square[W'`W`V'`V;`f`f'``]
 \efig$$
 and the bottom arrow is a weak isomorphism.
Clearly, $W'\to W$ is a bijection, so we need only show that
$\Hom(W,K)\to\Hom(W',K)$ is surjective.

I claim that $W'\inc W\x V'$ with the induced topology.  Let us define
$W''$ to be the subobject $W\x_VV'$ with the induced topology.  Since
$W'\to W$ and $W'\to V$ are continuous, the topology on $W'$ is at least
as fine as that of $W''$.  On the other hand, we do have $W''\to W$ and
$W''\to V'$ with the same map to $V$ so that we have $W''\to W'$, so
that the topology on $W''$ is at least as fine as that of $W'$.  Then we
have a commutative diagram
 $$\bfig
 \Square/>`into->`>`>/[W'`W`W\x V'`W\x V;``(\id,f)`]
 \efig$$
 Apply $\Hom(-,K)$ and use the injectivity of $K$ to get:
 $$\bfig
 \Square|aaaa|/<- `<<-`<-
`<- /[\Hom(W',K)`\Hom(W,K)`\Hom(W,K)\x\Hom(V',K)`
\Hom(W,K)\x\Hom(V,K);```\iso]
 \efig$$
 The bottom arrow is a bijection and the left hand arrow is a
surjection, which implies that the top arrow is a surjection.\epf

\prop $\sigma$ and $\tau$ are functors on $\Vsc$.\eth

\pf For $\sigma$, this is easy.  If $f:W\to V$ is a morphism, the
induced $\sigma f:\sigma W\to\sigma V$ will be continuous if and only if
its composite with every functional on $V$ is a functional on $W$, which
obviously holds.

To see that $\tau$ is a functor, suppose $f:W\to V$ is a morphism.  Form
the pullback
  $$\bfig
 \Square[W'`W`\tau V`V;`f'`f`]
 \efig$$
 It is a weak isomorphism by the preceding proposition.  Thus we get
$\tau W\to W'\to\tau V$.\epf

\prop If $V\to V'$ is a weak isomorphism, then $\sigma V\to \sigma V'$
and $\tau V\to\tau V'$ are isomorphisms.\eth

\pf For $\sigma$, this is obvious.  Clearly, $\tau V\to V\to \tau V'$ is
also a weak isomorphism so that $\tau V$ is one of the factors in the
computation of $\tau V'$ and then $\tau V'\to\tau V$ is a continuous
bijection, while the other direction is evident.\epf

\cor \label{WI}Both $\sigma$ and $\tau$ are idempotent, while
$\sigma\tau\iso\sigma$ and $\tau\sigma\iso\tau$.\epf\eth

\prop For any $V,V'\in\Vsc$, we have $\Hom(\sigma V,\sigma
V')\iso\Hom(\tau V,\tau V')$.\epf

\pf It is easiest to assume that the underlying objects $|V|=|\sigma
V|=|\tau V|$ and similarly for $V'$.  Then for any $f:V\to V'$, we also
have that $|f|=|\sigma f|=|\tau f|$.  Thus the two composition of the
two maps below
 $$\Hom(\sigma V,\sigma V')\to\Hom(\tau\sigma V,\tau\sigma
V')=\Hom(\tau V,\tau V')$$
 and
 $$\Hom(\tau V,\tau V')\to\Hom(\sigma \tau V,\sigma\tau V')\iso
\Hom(\sigma V,\sigma V')$$
 give the identity in each direction.\eth

Let $\Vsc_w\inc\Vsc$ and $\Vsc_s\inc\Vsc$ denote the full subcategories
of weak and strong objects, respectively. Then as an immediate corollary
to the preceding, we have:

\thm\label{weak=strong} $\tau:\Vsc_w\to\Vsc_s$ and
$\sigma:\Vsc_s\to\Vsc_w$ determine inverse equivalences of
categories.\epf\eth

\section{Chu and chu}\label{chu}

Now we add to the assumptions on $\Asc$ that it be a symmetric monoidal
closed category in which the underlying set of $A\r B$ is $\Hom(A,B)$.
We denote by \Esc and \Msc the classes of surjections and injections,
respectively.

We briefly review the categories $\Chu(\Asc,K)$ and $\chu(\Asc,K)$.
See \cite{B98} for details.  The
first has a objects pairs $(A,X)$ of objects of $\Asc$ equipped with a
``pairing'' $<-,->:A\* X\to K$.  A morphism $(f,g):(A,X)\to(B,Y)$
consists of a map $f:A\to B$ and a map $g:Y\to X$ such that
 $$\bfig
 \Square[A\*Y`B\*Y`A\*X`K;f\*Y`A\*g`<-,->`<-,->]
 \efig$$
 commutes.  This says that $<fa,y>=<a,gy>$ for all $a\in A$ and $y\in
Y$.  This can be enriched over $\Asc$ by internalizing this definition
as follows.  Note first that the map $A\*X\to K$ induces, by exponential
transpose, a map $X\to A\r K$.  This gives a map $Y\r X\to Y\r(A\r
K)\iso A\*Y\r K$.  There is a similarly defined arrow ${A\r B\to A\*Y\r
K}$. Define $[(A,X),(B,Y)]$ so that
 $$\bfig
 \Square[{[(A,X),(B,Y)]}`A\r B`Y\r X`A\*Y\r K;```]
 \efig$$
 is a pullback.  Then define
 $$(A,X)\r(B,Y)=([(A,X),(B,Y)],A\*Y)$$
 with $<(f,g),a\*y>=<fa,y>=<a,gy>$ and
 $$(A,X)\*(B,Y)=(A\*B,[(A,X),(Y,B)])$$
 with pairing $<a\*b,(f,g)>=<b,fa>=<a,gb>$.  The duality is given by
$(A,X)^*=(X,A)\iso(A,X)\r(K,\top)$
 where $\top$ is the tensor unit of $\Asc$.  Incidentally, the tensor
unit of $\Chu(\Asc,K)$ is $(\top,K)$.

The category $\Chu(\Asc,K)$ is complete (and, of course, cocomplete).
The limit of a diagram is calculated using the limit of the first
coordinate and the colimit of the second.  The full subcategory
$\chu(\Asc,K)\inc\Chu(\Asc,K)$ consists of those objects $(A,X)$ for
which the two transposes of $A\*X\to K$ are injective homomorphisms.
When $A\mon X\r K$, the pair is called separated and when $X\mon A\r K$,
it is called extensional.  In the general case, one must choose a
factorization system $(\Esc,\Msc)$ and assume that the arrows in \Esc
are epic and that \Msc is stable under $\r$, but here these conditions
are clear.  Let us denote by $\Chu_s(\Asc,K)$ the full subcategory of
separated pairs and by $\Chu_e(\Asc,K)$ the full subcategory of
extensional pairs.

The inclusion $\Chu_s(\Asc,K)\incto \Chu(\Asc,K)$ has a left adjoint $S$
and the inclusion $\Chu_e(\Asc,K)\incto\Chu(\Asc,K)$ has a right adjoint
$E$.  Moreover, $S$ takes an extensional pair into an extensional one
and $E$ does the dual.  In addition, when $(A,X)$ and $(B,Y)$ are
separated and extensional, $(A,X)\r(B,Y)$ is separated but not
necessarily extensional and, dually, $(A,X)\*(B,Y)$ is extensional, but
not necessarily separated.  Thus we must apply the reflector to the
internal hom and the coreflector to the tensor, but everything works out
and $\chu(\Asc,K)$ is also $*$-autonomous.  See \cite{B98} for details.

In the chu category, one sees immediately that in a map $(f,g):(A,X)\to
(B,Y)$, $f$ and $g$ determine each other uniquely.  So a map could just
as well be described as an $f:A\to B$ such that $x.\tilde y\in X$ for
every $y\in Y$.  Here $\tilde y:B\to K$ is the evaluation at $y\in Y$ of
the exponential transpose $Y\to B\r K$.

Although the situation in the category of abelian groups is as
described, in the case of vector spaces over a field, the hom and tensor
of two separated extensional pairs turns out to be separated and
extensional already (\cite{B96}).

\section{The main theorem}

\thm Assume the hypotheses of Theorem \ref{s&t} and also assume that
the canonical map $I\to K\r K$ is an isomorphism.  Then the categories
of weak spaces and strong spaces are equivalent to each other and to
chu$(\Asc,K)$ and are thus $*$-autonomous.\eth

\pf The first claim is just Theorem \ref{equiv}.  Now
define $F:\Vsc\to \chu$ by $F(V)=(|V|,\Hom(V,K))$ with evaluation as
pairing.  We first define the right adjoint $R$ of $F$.  Let $R(A,X)$ be
the object $A$, topologized as a subobject of $K^X$.  Since it is
already inside a power of $K$, it has the weak topology.  Let $f:|V|\to
A$ be a homomorphism such that for all $x\in X$,\quad$\tilde x.f\in
\Hom(V,K)$.  This just means that the composite $V\to R(A,X)\to
K^X\to^{\pi_x} K$ is continuous for all $x\in X$, exactly what is
required for the map into $R(A,X)$ to be continuous.  The uniqueness of
$f$ is clear and this establishes the right adjunction.

We next claim that $FR\iso\Id$.  That is equivalent to showing that
$\Hom(R(A,X),K)=X$.  Suppose $\phi:R(A,X)\to K$ is a functional.  By
injectivity, it extends to a $\psi: K^X\to K$.  It follows from
\ref{main}, there is a finite set of elements $x_1,\ldots,x_n\in X$ and
morphisms $\theta_1,\ldots,\theta_n$ such that $\psi$ factors as
$K^X\to K^n\to^{(\theta_1,\ldots,\theta_n)} K$.  Applied to $R(A,X)$,
this means that $\phi(a)=<\theta_1x_1,a>+\cdots+<\theta_nx_n,a>$.  But
the $\theta_i\in I$ and the tensor products are over $I$ so that the
pairing is a homomorphism $A\*_IX\to K$.  This means that
$\phi(a)=<\theta_1x_1+\cdots\theta_nx_n,a>$ and
$\theta_nx_1+\cdots+\theta_nx_n\in X$.

Finally, we claim that $RF=S$,
the left adjoint of the inclusion $\Vsc_w\inc \Vsc$.  If $V\in \Vsc$,
then $RFV=R(|V|,\Hom(V,K))$ which is just $V$ with the weak topology it
inherits from $K^{\Hom(V,K)}$, exactly the definition of $SV$.  It
follows that  $F|\Vsc_w$ is an equivalence.

Since $\Vsc_w$ and $\Vsc_s$ are equivalent to a $*$-autonomous category,
they are $*$-autonomous.\epf

The fact that the categories of weak and Mackey spaces are equivalent
was shown, for the case of B (Banach) spaces in \cite[Theorem 15, p.
422]{DS}. Presumably,
the general case has also been long known, but I am not aware of a
reference.

\section{Examples.}

\subsection*{Example 1. Vector spaces over a locally compact field.}
For the second example, let $K$ be a locally compact field.  Locally
compact fields have been classified, see \cite{Pont} or \cite{W}.
Besides the discrete fields and the real and complex numbers, they come
in two varieties.  The first are finite algebraic extensions of the
field $\Q_p$ which is the completion of the rational field in the
$p$-adic norm.  The second are finite algebraic extensions of the field
$\S_p$, which is the completion in $t$-adic norm of the field
$\Z_p\{t\}$ of Laurent series over the field of $p$ elements.  Notice
that all these locally compact fields are normed.

We take for $\Ssc$ the category of normed linear $K$-spaces, except in
the case that $K$ is discrete, we require also that the spaces have the
discrete norm.  We know that $K$ is injective in the discrete case.  The
injectivity of $K$ in the real or complex case is just the Hahn-Banach
theorem, which has been generalized ultrametric fields
according to the following, found in \cite{Rob}.

\thm[Ingelton]
Let $K$ be a spherically complete ultrametric field.  $E$ a $K$-normed
space and $v$ a subspace of $E$.  For every bounded linear functional
$\phi$ defined on $V$, there exists a bounded linear functional $\psi$
defined on $E$ whose restriction to $V$ is $\phi$ and such that
$\|\phi = \|\psi$.\eth

An ultrametric is a metric for which the ultratriangle inequality,
$\|{x+y}\le \|x \vee \|y$, holds.  This is obviously true for $p$-adic
and $t$-adic norms.  Spherically complete means that the meet of any
descending sequence of closed balls is non-empty.  This is known to be
satisfied by locally compact ultrametric spaces.

Regardless of the topology on a field $K$ (assuming it is topological
field), $K$ is its own endomorphism ring.

Notice that if $K$ is non-discrete, then what we have established is
that both $\Vsc_s$ and $\Vsc_w$ are equivalent to
$\chu(\Vect$-$|K|,|K|)$.  But exactly the same considerations show that
the same is true if we ignore the topology on $K$ and use the discrete
norm.  The category $\Ssc$ will now be the category of discrete
finite-dimensional $|K|$-vector spaces.  Its product and subobject
closure will consist of spaces that are mostly not discrete, but there
are still full subcategories of weakly and strongly topologized spaces
within this category and they are also equivalent to
$\chu(\Vect$-$|K|,|K|)$.

Thus, these categories really do not depend on the topologies.  Another
interpretation is that this demonstrates that, for these spaces, the
space of functionals replaces the topology, which was arguably Mackey's
original intention.

\subsection*{Example 2. Locally compact abelian groups.}
For the abelian groups, we take for $\Vsc$ the category of those abelian
that a subgroups (with the induced topology) of products of locally
compact abelian groups.  The object $K$ in this case is the circle group
$\R/\Z$.  A simple representation of this group is as the closed
interval $[-1/2,1/2]$ with the endpoints identified and addition mod 1.
The group is compact.  Let $U$ be the open interval $(-1/3,1/3)$.  It is
easy to see that any non-zero point in that interval, added to itself
sufficiently often, eventually escapes that neighborhood so that $U$
contains no non-zero subgroup.  It is well-known that the endomorphism
group of the circle is $\Z$.

If $f:G\to K$ is a homomorphism such that $T=f\inv(U)$ is open in $G$,
let $T=T_1,T_2,\ldots,T_n,\ldots$ be a sequence of open sets in $G$ such
that $T_{i+1}+T_{i+1}\inc T_i$ for all $i$.  Let
$U_i=(-2^{-i}/3,2^{-i}/3)\inc K$.  Then the $\{U_i\}$ form a
neighborhood base in $K$ and one readily sees that $f\inv(U_i)\inc T_i$
which implies that $f$ is continuous.

We take for $\Ssc$ the category of locally compact abelian groups.  The
fact that $K$ is an injective follows directly from the
Pontrjagin duality theorem.  A result \cite[Theorem 1.1]{G} says that
every locally compact group is strongly topologized.  Thus both
categories of weakly topologized and strongly topologized groups that
are subobjects of products of locally compact abelian groups are
equivalent to $\chu(\Ab,|K|)$ and thus are *-autonomous.

We can ask if the same trick of replacing $K=\R/\Z$ by $|K|$, as in the
first example, can work.  It doesn't appear so.  While $\Hom(K,K)=\Z$,
the endomorphism ring of $|K|$ has cardinality $2^c$ and is
non-commutative, so we cannot draw no useful inference about maps from
$|K|^n\to |K|$, even for finite $n$.

\subsection*{Example 3. Modules over a self injective cogenerator.}
 If we examine the considerations that are used in vector spaces over a
field, it is clear that what is used is that a field is both an
injective module over itself and a cogenerator in the category of
vector spaces.  Then if $K$ is a such a commutative ring, we can let
\Tsc be the category of topological $K$-modules, \Ssc be the full
subcategory of submodules of finite powers of $K$ with the discrete
topology and \Vsc the limit closure of \Ssc.  Then $\chu(\Mod_K,K)$ is
equivalent to each of the categories $\Vsc_s$ and $\Vsc_w$ of
topological $K$-modules that are strongly and weakly topologized,
respectively, with respect to their continuous linear functionals into
$K$.

We now show that there is a class of commutative rings with that
property.  Let $k$ be a field and $K=k[x]/(x^n)$.  When $n=2$, this is
called the ring of dual numbers over $k$.

\prop \label{sinj}$K$ is self injective.\eth

We base this proof on the following well-known fact:

\lem Let $k$ be a commutative ring, $K$ is a $k$-algebra, $Q$ an
injective $k$-module, and $P$ a flat right $K$-module then $\Hom_k(K,Q)$
is an injective $K$-module.\eth

The $K$-module structure on the Hom set is given by $(rf)(a)=f(ar)$ for
$r\in K$ and $a\in P$.

\pf Suppose $A\mon B$ is an injective homomorphism of $K$-modules.  Then
we have
 $$\bfig
\square(0,500)/>`>`>` >->/<1500,500>[\Hom_R(B,\Hom_k(P,Q))`\Hom_R(A,
\Hom_k(P,Q))` \Hom_k(P\*_RB,Q)`\Hom_k(P\*_RA,Q);`\iso`\iso`]
 \efig$$
 and the flatness of $P$, combined with the injectivity of $Q$ force the
bottom arrow to be a surjection.\epf

\pf[of \ref{sinj}]  From the lemma it follows that $\Hom_k(K,k)$ is a
$K$-injective.  We claim that, as $K$-modules, $\Hom_k(K,k)\iso K$.  To
see this, we map $f:K\to\Hom_k(K,k)$.  Since these are vector spaces
over $k$, we begin with a $k$-linear map and show it is $K$-linear.  A
$k$-basis for $K$ is given by $1,x,\ldots,x^{n-1}$.  We define
$f(x^i):K\to k$ for $0\le \le n-1$ by $f(x^i)(x^j)=\delta_{i+j,n}$ (the
Kronecker $\delta$).  For this to be $K$-linear, we must show that
$f(xx^i)=xf(x^i)$.  But
$$f(xx^i)(x^j)=f(x^{i+1}(x^j))=\delta_{i+1+j,n}=f(x^i)(x^{j+1})=
(xf(x^i))(x^j)$$
 Clearly, the $f(x^i)$, for $0\le i\le n$ are linearly independent and
so $f$ is an isomorphism.\epf

\prop $K$ is a cogenerator in the category of $K$-modules.\eth

\pf Using the injectivity, it suffices to show that every cyclic module
can be embedded into $K$.  Suppose $M$ is a cyclic module with generator
$m$.  Let $i$ be the first power for which $x^im=0$.  I claim that
$m,xm,\ldots,x^{i-1}m$ are linearly independent over $k$.  If not,
suppose that $\lam_0m+\lam_1xm+\cdots\lam_{i-1}x^{i-1}m=0$ and not all
coefficients 0.  Let $\lam_j$ be the first non-zero coefficient, so that
$\lam_jx^j+\cdots+\lam_{i-1}x^{i-1}m=0$.  Multiply this by $x^{i-j-1}$
and use that $x^lm=0$ for $l\ge i$ to get $\lam_jx^{i-1}m=0$.  But by
assumption, $x^{i-1}m\neq0$ so that this would imply that $\lam_j=0$,
contrary to hypothesis.  Thus there is a $k$-linear map $f:M\to K$ given
by $f(x^jm)=x^{n-i+j}$.  Since the $x^j$ are linearly independent, this
is $k$-linear and then it is clearly $K$-linear.\epf

 \comment
\subsection*{Modules over $k[x]/(x^2)$} Let $k$ be a field and
$K=k[x]/(x^2)$ the so-called ring of dual numbers over $k$.  It is
well-known that $K$ is $K$-injective.  We claim it is also cogenerates
the category of $K$

\section{Concluding remarks}

There is a curious conclusion to be drawn from all this.  Consider, for
example, the case of real vector spaces.  We could treat this using for
$\Ssc$ the category of normed vector spaces.  Or we could use the
category of finite dimensional vector spaces with the usual Euclidean
topology.  Or we could just treat $\R$ as an abstract field with the
discrete topology.  For each choice we get strong and weak topologies.
This gives us six distinct $*$-autonomous categories (at least, there
may be other possibilities).  But all of them are equivalent, because
they are all equivalent to $\chu(\Vsc,\R)$, where $\Vsc$ is the category
of real vector spaces.  In other words, the topology of the reals and of
various vector spaces plays no role in the structure of these
categories, qua categories.  The same thing happens with the complex
numbers and with the other locally compact fields.  There are also some
choices with topological abelian groups (see \cite{BK}) that all lead to
$\chu(\Ab,K)$.

Pairs equipped with a bilinear map were originally introduced by Mackey
\cite{M}.  Although his motivation is not entirely clear, it is at least
plausible that in a pair $(E,E')$ he thought of $E'$ not so much as a
set of continuous linear functionals with respect to a topology on $E$
(although he did think of them that way) but perhaps even more as
embodying a replacement for the very idea of a topology.  However, he
appears never to have defined what he would mean by a continuous
homomorphism in this setting, although it is clear enough what it would
have to be.  Even less did he suggest any thought of an space of such
transformations between two such spaces  or of a tensor product of them.

\endcomment

\section{Interpretation of the dual of an internal hom}
 These remarks are especially relevant to the vector spaces, although
they are appropriate to the other examples.  The fact that $(U\r V)^*
\iso U \r V^*$ can be interpreted that the dual of $U\r V$ is a subspace
of $V\r U$, namely those linear transformations of finite rank.  An
element of the form $u\*v^*$ acts as a linear transformation by the
formula $(u\*v^*)(v)=<v,v^*>u$.  This is a transformation of row rank 1.
Sums of these elements is similarly an element of finite rank.

This observation generalizes the fact that in the category of finite
dimensional vector spaces, we have that $(U\r V)^* \iso V\r U$ (such a
category is called a compact $*$-autonomous category).  In fact, Halmos
avoids the complications of the definition of tensor products in that
case by \emph{defining} $U\* V$ as the dual of the space of
bilinear forms on $U\+ V$, which is quite clearly
equivalent to the dual of $U\r V^*\iso V\r U^*$ (\cite[Page
40]{H}).  (Incidentally, it might be somewhat pedantic to point out
that Halmos's definition makes no sense since $U\+ V$ is a vector
space in its own right and a bilinear form on a vector space is absurd.
It would have been better to use the equivalent form above or to define
Bilin$(U,V)$.)

Since linear transformations of finite rank are probably not of much
interest in the theory of topological vector spaces, this may explain
why the internal hom was not pursued.

\section{Appendix: Some generalities on adjoints.}

In an earlier version of this paper, the proof of \ref{weak=strong} was
based on some formal results of adjoints.  The argument got greatly
simplified and these results were not needed, largely because of the
concreteness of the categories involved.  Still, it seemed worthwhile to
include these formal results.

\prop \label{2adj}Suppose $F:\Asc\to\Bsc$ is a functor that has both a
left adjoint $L$ and a right adjoint $R$.  Then $L$ is full and faithful
if and only $R$ is.\eth

\pf Suppose that $L$ is full and faithful.  Then we have, for any
$B,B'\in\Bsc$,
 $$\Hom(B,B')\iso\Hom(LB,LB')\iso\Hom(B,FLB')$$
 so that, by the Yoneda lemma, the front adjunction $B'\to FLB'$ is an
isomorphism.  Then
 $$\Hom(B,B')\iso\Hom(FLB,B')\iso\Hom(LB,RB')\iso\Hom(B,FRB')$$
 which implies that the back adjunction $FRB'\to B'$ is also an
isomorphism, which is possible only if $R$ is full and faithful.  The
reverse implication is just the dual.\epf

\prop Suppose
 $$\bfig
 \Square[A`B`C`D;f`g`h`k]
 \efig$$
 is a commutative square and both $f$ and $k$ are isomorphisms.  Then
 $$\bfig
 \Square/<-`>`>`<-/[A`B`C`D;f\inv`g`h`k\inv]
 \efig$$
 also commutes.  If, in addition, $g$ and $h$ are isomorphisms, then
 $$\bfig
 \Square/<-`<-`<-`<-/[A`B`C`D;f\inv`g\inv`h\inv`k\inv]
 \efig$$
 also commutes.\eth

\pf From $gf=kh$, we infer that $k\inv gff\inv=k\inv khf\inv$, which is
the first assertion and the second is immediate.\epf

\thm Suppose \Csc is a category, $I:\Bsc\to\Csc$ the inclusion of a full
subcategory with a left adjoint $S$ and $J:\Dsc\to\Csc$ is the inclusion
of a full subcategory with a right adjoint $T$.  Let $\alpha:1\to TS$
and $\beta:SI\to^{\iso}1$ be the front and back adjunctions for $S\adj
I$ and $\delta:1\to^{\iso}TJ$ and $\epsilon:JT\to 1$ do the same for
$J\adj T$.  Suppose, in addition, that $IS\eps:ISJT\to IS$ and
$JT\alpha:JT\to JTIS$ are isomorphisms.  Then $JT\adj IS$.\eth

\pf If $f:JT C\to C'$ is given, define $\mu f:C\to ISC'$
as the composite
 $$C\to^{\alpha C}ISC\to^{(IS\eps C)\inv}\to^{ISf}ISC'$$
 If $g:C\to ISC'$ is given, define $\nu g:JTC\to C' $ as the composite
 $$JTC\to^{JTg}JTSIC'\to^{(JT\alpha C')\inv}JTC'\to^{\eps C'}C'$$
 We must show that $\mu$ and $\nu$ are inverse operations.  The upper
and right hand arrows calculate $\nu\mu f$ and the squares commute by
naturality of by applying the preceding proposition to a naturally
commuting square.
 $$\bfig
 \node 1a(0,1000)[JTC]
 \node 1b(700,1000)[JTISC]
 \node 1c(1600,1000)[JTISJTC]
 \node 1d(2500,1000)[JTISC']
 \node 2b(700,500)[ISC]
 \node 2c(1600,500)[ISJTC]
 \node 2d(2500,500)[ISC']
 \node 3c(1600,0)[JTC']
 \node 3d(2500,0)[C']
 \arrow[1a`1b;JT\alpha C]
 \arrow[1b`1c;JT(IS\eps)\inv C]
 \arrow[1c`1d;JTISf]
 \arrow|l|[1a`2b;\iso]
 \arrow|r|[1b`2b;(JT\alpha)\inv C]
 \arrow|r|[1c`2c;(JT\alpha)\inv JTC]
 \arrow|r|[1d`2d;JT\alpha C']
 \arrow[2c`2d;ISf]
 \arrow[2b`2c;(IS\eps)\inv C]
 \arrow|l|[2b`3c;\iso]
 \arrow|r|[2c`3c;IS\eps C']
 \arrow|r|[2d`3d;\eps C']
 \arrow[3c`3d;f]
 \efig$$
 This shows that $\nu\mu f=f$ and $\mu\nu g=g$ is handled similarly.
Thus $\Hom(JTC,C')\iso\Hom(C,ISC')$

\prop $JTI\adj S$ and $T\adj ISJ$.\eth

\pf We have $\Hom(JTIB,C)\iso\Hom(IB,ISC)\iso\Hom(B,SC)$ since $I$ is
full and faithful.  The second one is proved similarly.\epf

\prop $JTI:\Bsc\to\Csc$ and $ISJ:\Dsc\to\Bsc$ are full and faithful.\eth
 \pf $S$ has a right adjoint $I$ and a right adjoint $JTI$.  Since $I$
is full and faithful, so is $JTI$ be \ref{2adj}.  The argument for $ISJ$
is similar.\epf

\cor $TI:\Bsc\to \Dsc$ is left adjoint to $SJ:\Dsc\to \Csc$ and each is
an equivalence.\eth

\pf $\Hom(TIB,D)\iso\Hom(JTIB,JD)\iso\Hom(B,SJD)$, gives the adjunction.
Moreover, $\Hom(TIB,TIB')\iso\Hom(JTIB,JTIB')\iso\Hom(B,B')$ since $JTI$
is full and faithful and a similar argument works for $SJ$.\epf

Applying to the results of Section \ref{weak=strong}, we conclude that:

\thm \label{equiv}The functors $TI:\Vsc_w\to\Vsc_s$ and
$SJ:\Vsc_s\to\Vsc_w$ are adjoint equivalences.\epf\eth

\subsection{Application}
 This was originally applied to the proof of \ref{weakstrong} as
follows.

Let $I:\Vsc_s\to \Vsc$ and $J:\Vsc_w\to\Vsc$ denote the inclusions of
into $\Vsc$ of the full subcategories
consisting of the weak and strong objects, respectively.

\thm The functor $S:\Vsc\to \Vsc_w$ for which $SV=\sigma V$ is left
adjoint to $I$.  Similarly, the functor $T:\Vsc\to\Vsc_s$ for which
$TV=\tau V$ is right adjoint to $J$.\eth

\pf First we note that for $\sigma I=ISI\iso I$.  Then for any
$V\in\Vsc$ and $V'\in\Vsc_w$ we have the composite
 $$\Hom(V,IV')\to\Hom(\sigma V,\sigma IV')\to\Hom(V,IV')$$
 is the identity and the second arrow is an injection, so that both
arrows are isomorphisms.  Thus $\Hom(V,IV)\iso\Hom(\sigma V,IV)$.  Then
we have
 $$\Hom(V,IV)\iso\Hom(\sigma V,IV')=\Hom(ISV,IV')\iso\Hom(SV,V')$$
 The second assertion is dual.\epf

\cor\label{WI'} $SITJ\to SI$ using the back adjunction and $TJ\to TJSI$
using the front adjunction are isomorphisms.\eth

\pf This is immediate from Corollary \ref{WI}.

\begin{references*}

\bibitem[Arens, 1947]{A} Arens, R. (1947), \emph{Duality in linear
spaces}.  Duke Math. J. {\bf14}, 787--794.  M.R. 0022651.

\bibitem [Barr, 1996]{B96} M.\ Barr (1996), \emph{Separability of tensor
in Chu categories of vector spaces}.  \mscs6, 213--217.

\bibitem [Barr, 1998]{B98} M.\ Barr (1998), \emph{The separated
extensional Chu category}.  \tac{4}, 127--137.

\bibitem [Barr, 2000]{B00} M.\ Barr (2000), \emph{On $*$-autonomous
categories of topological vector spaces}.  \cahiers {41}, 243--254.

\bibitem [Barr \& Kleisli, 2001]{BK} M.\ Barr and H.\ Kleisli (2001),
\emph{On Mackey topologies in topological abelian groups}.  \tac{9},
54-62.

\bibitem[Barr, 2006]{top} M. Barr (2006), Topological $*$-autonomous
categories.  \tac{16} (2006), 700--708.

\bibitem [Dunford \& Schwartz, 1958]{DS} N. Dunford and J. Schwartz
(1958), Linear Operators, Interscience.

\bibitem[Glicksberg, 1962]{G} I.\ Glicksberg (1962), \emph{Uniform
boundedness for groups}.  Canad.  Math.  J. {\bf 14}, 269--276.

\bibitem[Halmos, 1958]{H} P. Halmos (1958), \emph{Finite Dimensional
Vector Spaces}, second edition. D. Van Nostrand.

\bibitem[Mackey, 1945]{M} G.\ Mackey (1945), \emph{On infinite
dimensional vector spaces}.  Trans.\ Amer.\ Math.\ Soc.\ {\bf57},
155--207.

\bibitem[Mackey, 1946]{M46} G.\ Mackey (1946), \emph{On convext
topological linear spaces}.  Trans. Amer. Math. Soc. {\bf 66}, 519--537.
M.R. 0020214.

\bibitem[Pontrjagin, 1968]{Pont}
Pontryagin, Lev S. (1986).  Topological Groups. trans. from Russian by
Arlen Brown and P.S.V.  Naidu (3rd ed.).  New York:  Gordon and Breach
Science Publishers.

\bibitem[Robert, 2000]{Rob} Alain M. Robert (2000),  A course in p-adic
analysis.  Graduate texts in Math. {\bf 198}, Springer-Verlag.

\bibitem[Weil, 1967]{W} Andr\'e Weil (1967), Basic Number Theory.
Grundlehren {\bf 144}, Springer-Verlag.

\end{references*}
\end{document}